\pdfoutput=1
\documentclass[11pt,a4paper,american]{article}
\usepackage{babel}
\usepackage[T1]{fontenc}
\usepackage[utf8]{inputenc}
\usepackage{fullpage}
\usepackage[babel]{microtype}
\usepackage{csquotes}
\usepackage{authblk}
\usepackage{titlesec}
\usepackage[giveninits,sortcites,maxbibnames=10,style=ext-numeric]{biblatex}
\usepackage{hyperref}
\usepackage{enumitem}
\usepackage{amsmath}
\usepackage{mathtools}
\usepackage{amsthm}
\usepackage{thmtools}
\usepackage{thm-restate}
\usepackage{dsfont}
\usepackage{amssymb}
\usepackage{xcolor}
\usepackage[noabbrev,capitalise,sort]{cleveref}
\usepackage{mleftright}

\hypersetup{pdfencoding=auto,colorlinks=true,linkcolor=blue,citecolor=teal,hypertexnames=false,bookmarksnumbered=true,pdfauthor={Marco Pozza},pdftitle={Aubry Set of Eikonal Hamilton–Jacobi Equations on Networks},pdfkeywords={Hamilton–Jacobi equations;Aubry set;embedded networks;comparison principle;Hopf–Lax formula},pdfsubject={MSC(2020): 35F21, 35R02, 35B51, 49L25}}

\titlespacing{\paragraph}{0pt}{1em}{.7em}

\declaretheorem[style=plain,parent=section,title=Theorem,refname={Theorem,Theorems}]{theo}
\declaretheorem[style=plain,sibling=theo,title=Proposition,refname={Proposition,Propositions}]{prop}
\declaretheorem[style=plain,sibling=theo,title=Corollary,refname={Corollary,Corollaries}]{cor}
\declaretheorem[style=plain,sibling=theo,title=Lemma,refname={Lemma,Lemmas}]{lem}
\declaretheorem[style=definition,sibling=theo,title=Definition,refname={Definition,Definitions}]{defin}
\declaretheorem[style=remark,sibling=theo,title=Remark,refname={Remark,Remarks}]{rem}

\AddToHook{cmd/appendix/before}{\crefalias{chapter}{appendix}\crefalias{section}{appendix}\crefalias{subsection}{appendix}}

\newlist{Hassum}{enumerate}{1}
\setlist[Hassum]{label=\textbf{(H\arabic*)},ref=\textnormal{\textbf{(H\arabic*)}},font=\normalfont}
\newlist{myenum}{enumerate}{1}
\setlist[myenum]{label=\textbf{\roman*)},ref=\textnormal{\textbf{(\roman*)}},font=\normalfont}
\newlist{mylist}{itemize}{1}
\setlist[mylist]{label=\textbullet,font=\normalfont}

\crefname{equation}{}{}
\crefname{Hassumi}{condition}{conditions}
\crefname{myenumi}{item}{items}

\DeclareFieldFormat{titlecase:title}{\MakeSentenceCase*{#1}}

\DeclareMathOperator{\interior}{int}

\newcommand{\Acal}{\mathcal A}
\newcommand{\Ccal}{\mathcal C}
\newcommand{\Ecal}{\mathcal E}
\newcommand{\Hcal}{\mathcal H}
\newcommand{\Nds}{\mathds N}
\newcommand{\Ocal}{\mathcal O}
\newcommand{\ovOcal}{\overline\Ocal}
\newcommand{\Rds}{\mathds R}
\newcommand{\Vbf}{\mathbf V}
\newcommand{\ovs}{\overline s}
\newcommand{\ovx}{\overline x}
\newcommand{\wtgamma}{\widetilde\gamma}
\newcommand{\ovmu}{\overline\mu}
\newcommand{\myemail}[2][]{\textsuperscript{#1}\href{mailto:#2}{\texttt{#2}}}
\newcommand{\hrefc}{\textup{(}\hyperref[eq:globeik]{\ensuremath{\Hcal}\textup{J}\ensuremath{c}}\textup{)}}
\newcommand{\ghrefc}{\textup{(}\hyperref[eq:eikg]{\textup{HJ}\textsubscript{\ensuremath{\gamma}}\ensuremath{c}}\textup{)}}

\newcommand{\abscite}[2][]{\citeauthor{#2} \parentext{\citefield{#2}[journaltitle]{shortjournal}, \bibhyperref[#2]{\citeyear{#2}}}}

\title{Aubry Set of Eikonal Hamilton--Jacobi Equations on Networks}
\author{Marco Pozza}
\affil{Link Campus University, Rome, Italy. \textit{Email address:} \myemail{m.pozza@unilink.it}}
\date{}

\begin{document}

    \maketitle

    \begin{abstract}
        We extend the study of eikonal Hamilton--Jacobi equations posed on networks performed by \abscite{SiconolfiSorrentino18} to a more general setting. Their approach essentially exploits that such equations correspond to discrete problems on an abstract underlying graph. However, a specific condition they assume can be rather restricting in some settings, which motivates the generalization we propose. We still get an \emph{Aubry set}, which plays the role of a uniqueness set for our problem and appears in the representation of solutions. Exploiting it we establish a new comparison principle between super and subsolutions to the equation.
    \end{abstract}

    \paragraph{2020 Mathematics Subject Classification:} 35F21, 35R02, 35B51, 49L25.

    \paragraph{Keywords:} Hamilton--Jacobi equations, Aubry set, embedded networks, comparison principle, Hopf--Lax formula.

    \section{Introduction}

    The aim of this paper is to extend the study of eikonal-type Hamilton--Jacobi equations on networks performed in~\cite{SiconolfiSorrentino18} to a more general setting. We establish a comparison principle and provide a Hopf--Lax type representation formula for solutions of such problems.

    Since the pioneering works of~\cite{CamilliSchieborn12}, there is an increasing interest in the study of Hamilton--Jacobi equations posed on networks, see for instance~\cite{ImbertMonneau17,Morfe20,LionsSouganidis17,CamilliMarchi23,AchdouTchou15,AchdouCamilliCutriTchou12,ImbertMonneauZidani12}. This attention is motivated by several applications (traffic models, data transmission, management of computer cluster, etc.) as well as the analysis of some theoretical issues related to the discontinuities of the Hamiltonian in this framework, as pointed out in the recent comprehensive monograph~\cite{BarlesChasseigne24}. Eikonal equations are particularly important due to their connections with other problems. For example, solutions to discounted and time-dependent Hamilton--Jacobi equations converge, under suitable conditions, to solutions to corresponding eikonal equations as the discount factor goes to 0 (\cite{PozzaSiconolfi21}) and the time positively diverges (\cite{Pozza23}), respectively.

    We consider a connected network $\Gamma$ embedded in $\Rds^N$ with a finite number of \emph{arcs} $\gamma$, namely regular curves parametrized in $[0,1]$, linking points of $\Rds^N$ called \emph{vertices}. A Hamiltonian on $\Gamma$ is a collection of Hamiltonians $H_\gamma:[0,1]\times\Rds\to\Rds$ indexed by the arcs, with the crucial feature that Hamiltonians associated to arcs possessing different support are totally unrelated. We are interested in the corresponding family of equations
    \begin{equation}\label{eq:protoeik}
        H_\gamma(s,Du)=a,\qquad \text{on }[0,1],
    \end{equation}
    where $a$ is a constant independent of $\gamma$. A solution to such problem on $\Gamma$ is a continuous function $u:\Gamma\to\Rds$ satisfying a suitable condition on the vertices and such that $u\circ\gamma$ is a solution, in the viscosity sense, to~\cref{eq:protoeik} for each arc $\gamma$. More specifically, a solution must satisfy at each vertex a state-constraint type boundary condition, see \cref{defsol}\ref{stateconst}. Such condition has also been considered for other stationary equations on networks (\cite{LionsSouganidis16,PozzaSiconolfi21}) and even time-dependent problems (\cite{Siconolfi22,ImbertMonneau17}), which are more involved and require additional conditions, like \emph{flux limiters}---namely upper bounds of the time derivative of subsolutions at vertices.

    We show that there is a unique value of the parameter $a$ in~\cref{eq:protoeik}, called \emph{critical value}, such that the eikonal problem admits solutions on the whole network. The critical value is in addition deeply related to the Aubry--Mather theory and the homogenization problem on networks, see~\cite{PozzaSiconolfiSorrentino24,SiconolfiSorrentino23}.

    Note that if, for every arc $\gamma$, $a$ is greater than
    \begin{equation*}
        a_\gamma\coloneqq\max\limits_{s\in[0,1]}\min\limits_{\mu\in\Rds}H_\gamma(s,\mu),
    \end{equation*}
    any local equation~\cref{eq:protoeik} admits unique solutions once the boundary values are fixed without any restriction, see \cite{SiconolfiSorrentino18}. This holds true also in the case where the arc is a \emph{loop}, i.e., a closed curve. This crucial local uniqueness property actually allows establishing, as done in~\cite{SiconolfiSorrentino18}, a correspondence of the equation on the network with a discrete problem on an underlying graph. When instead $a=a_\gamma$, it can happen that the local problem indexed by $\gamma$ with fixed boundary values admits multiple solutions. In this frame, to keep the correspondence with the discrete problem, condition~\ref{condmina} is introduced in~\cite{SiconolfiSorrentino18}. This assumption---namely the request that, for each $\gamma$ with $a_\gamma$ equal to the critical value, the minimum of $H_\gamma$ is equal to $a_\gamma$---ensures that there is uniqueness set, called \emph{Aubry set}, which, according to the aforementioned discrete approach, roughly corresponds to a subset of the underlying graph. Once an admissible trace is assigned on such set, a unique critical solution is identified via a Hopf--Lax type representation formula.

    Condition~\ref{condmina} is usually not assumed in other works about Hamilton--Jacobi equations on networks. Even the ones cited before which establish some connection with the eikonal equation, impose~\ref{condmina} with the sole objective of proving said connection. This condition can also be rather restrictive in some applications, e.g., it is impractical for problems where the critical value is unknown beforehand.

    The main novelty of our work is that we overcome this restriction by extending the weak KAM analysis carried out in~\cite{SiconolfiSorrentino18} to a setting where condition~\ref{condmina} is not assumed, losing the correspondence with the discrete problem. We will exclusively work on the network, and to do so we need to adapt the tools developed in~\cite{SiconolfiSorrentino18} to a continuous setting. Not surprisingly, while we can still retrieve a uniqueness set, it does not correspond in general to a subset of the abstract graph. For example, it could be made up of a countable number of disjointed pieces of arcs or just a single point. Using this new Aubry set, we provide a Hopf--Lax representation formula for solutions and a comparison principle between super and subsolutions, which, to our knowledge, is a novelty for this problem.

    The article is organized as follows. We introduce the notions of network and Hamiltonian defined on it, together with our assumptions, in \cref{netsec}. In \cref{eiksec} we present the problem under investigation. \Cref{aubsec} is devoted to the characterization of the Aubry set and the representation formula for the solution. In \cref{compsec} we provide our comparison principle. \Cref{optcurve} contains the proofs of some auxiliary results.

    \paragraph{Acknowledgments.} The author is a member of the INdAM research group GNAMPA.

    \section{Networks}\label{netsec}

    Here we describe our framework, namely what is a network, a Hamiltonian defined on it and some related concepts useful for our study.

    \subsection{Basic Definitions}

    We fix a dimension $N$ and $\Rds^N$ as ambient space. An \emph{embedded network}, or \emph{continuous graph}, is a subset $\Gamma\subset\Rds^N$ of the form
    \begin{equation*}
        \Gamma=\bigcup_{\gamma\in\Ecal}\gamma([0,1])\subset\Rds^N,
    \end{equation*}
    where $\Ecal$ is a finite collection of regular (i.e., $C^1$ with non-vanishing derivative) simple oriented curves, called \emph{arcs} of the network, that we assume, without any loss of generality, parameterized on $[0,1]$. Note that we are also assuming existence of one-sided derivatives at the endpoints 0 and 1. We stress out that a regular change of parameters does not affect our results.

    Observe that on the support of any arc $\gamma$, we also consider the inverse parametrization defined as
    \begin{equation*}
        \widetilde\gamma(s)\coloneqq\gamma(1-s),\qquad \text{for }s\in[0,1].
    \end{equation*}
    We call $\widetilde\gamma$ the \emph{inverse arc} of $\gamma$. We assume that if $\gamma\in\Ecal$ then $\wtgamma\in\Ecal$ and
    \begin{equation}\label{eq:nosovrap}
        \gamma((0,1))\cap\gamma'([0,1])=\emptyset,\qquad \text{whenever }\gamma'\ne\gamma,\widetilde\gamma.
    \end{equation}

    We call \emph{vertices} the initial and terminal points of the arcs, and denote by $\Vbf$ the sets of all such vertices. It follows from~\cref{eq:nosovrap} that vertices are the only points where arcs with different support intersect and, in particular,
    \begin{equation*}
        \gamma((0,1))\cap\Vbf=\emptyset,\qquad \text{for any }\gamma\in\Ecal.
    \end{equation*}

    We assume that the network is connected, namely given two vertices there is a finite concatenation of arcs linking them.

    For each $x\in\Vbf$, we define $\Gamma_x\coloneqq\{\gamma\in\Ecal:\gamma(1)=x\}$.

    The network $\Gamma$ inherits a geodesic distance, denoted with $d_\Gamma$, from the Euclidean metric of $\Rds^N$. It is clear that given $x$, $y$ in $\Gamma$ there is at least a geodesic linking them. The geodesic distance is in addition equivalent to the Euclidean one.

    We also consider a differential structure on the network by defining, for every $x\in\Gamma$, the \emph{tangent space} to $\Gamma$ at $x$, $T_x\Gamma$ in symbols, as the set made up by the $q\in\Rds^N$ of the form
    \begin{equation*}
        q=\lambda\dot\gamma(s),\qquad \text{if $x=\gamma(s)$, $s\in[0,1]$, with $\lambda\in\Rds$}.
    \end{equation*}
    Note that $\dot\gamma(s)$ is univocally determined, up to a sign, if $x\in\Gamma\setminus\Vbf$ or in other words if $s\notin\{0,1\}$. We further define the \emph{cotangent space} $T^*_x\Gamma$ as the dual space of $T_x\Gamma$, and set the \emph{tangent bundle} $T\Gamma$ and \emph{cotangent bundle} $T^*\Gamma$ as the disjoint union of the $T_x\Gamma$ and $T_x^*\Gamma$ at all points of $\Gamma$, respectively.

    By curve we mean throughout the paper an \emph{absolutely continuous} curve. We point out that the pair $\mleft(\xi,\dot\xi\mright)$, where $\xi$ is a curve in $\Gamma$, is naturally contained in $T\Gamma$. Let $\xi:[0,T]\to\Gamma$ and $\xi':[0,T']\to\Gamma$ be two curves such that $\xi(T)=\xi'(0)$. We define their \emph{concatenation} as the curve $\xi*\xi':[0,T+T']\to\Gamma$ such that
    \begin{equation*}
        \xi*\xi'(t)\coloneqq\mleft\{
        \begin{aligned}
            &\xi(t),&& \text{if }t\in[0,T),\\
            &\xi'(t-T),&& \text{if }t\in\mleft[T,T+T'\mright].
        \end{aligned}
        \mright.
    \end{equation*}
    Notice that $*$ is an associative operation.

    \begin{defin}\label{diffdef}
        Given a function $f:\Gamma\to\Rds$ such that $f\circ\gamma$ is differentiable at $x=\gamma(s_0)$, where $\gamma\in\Ecal$ and $s\in(0,1)$, we define its differential at $x$, denoted by $D_\Gamma f(x)$, as the unique covector in $T^*\Gamma$ such that
        \begin{equation*}
            D_\Gamma f(x)\dot\gamma(s_0)\coloneqq\frac d{ds}f(\gamma(s))|_{s=s_0}.
        \end{equation*}
    \end{defin}

    \subsection{Hamiltonians on \texorpdfstring{$\Gamma$}{Γ}}

    A Hamiltonian on $\Gamma$ is a collection of Hamiltonians $\Hcal\coloneqq\{H_\gamma\}_{\gamma\in\Ecal}$, where
    \begin{alignat*}{2}
        H_\gamma:[0,1]\times\Rds&&\:\longrightarrow\:&\Rds\\
        (s,\mu)&&\:\longmapsto\:&H_\gamma(s,\mu),
    \end{alignat*}
    satisfying
    \begin{equation}\label{eq:condcompH}
        H_{\wtgamma}(s,\mu)=H_\gamma(1-s,-\mu),\qquad \text{for any arc }\gamma.
    \end{equation}
    We emphasize that, apart the above compatibility condition, the Hamiltonians $H_\gamma$ are \emph{unrelated}.

    We require any $H_\gamma$ to be:
    \begin{Hassum}
        \item continuous in both arguments;
        \item \label{condcoerc} coercive in $\mu$;
        \item \label{condqconv} strictly quasiconvex in $\mu$, which means that, for any $s\in[0,1]$, $\mu,\mu'\in\Rds$ and $\rho\in(0,1)$,
        \begin{equation*}
            H_\gamma\mleft(s,\rho\mu+(1-\rho)\mu'\mright)<\max\mleft\{H_\gamma(s,\mu),H_\gamma\mleft(s,\mu'\mright)\mright\}.
        \end{equation*}
        When \cref{condcoerc} holds the above assumption is equivalent to require, for every $a\in\Rds$ and $s\in[0,1]$, the convexity of the sublevel set $\{\mu\in\Rds:H_\gamma(s,\mu)\le a\}$ (provided it is nonempty) and
        \begin{equation}\label{eq:intsublevel}
            \interior\{\mu\in\Rds:H_\gamma(s,\mu)\le a\}=\{\mu\in\Rds:H_\gamma(s,\mu)<a\},
        \end{equation}
        where $\interior$ denotes the interior of a set.
    \end{Hassum}

    We define the support functions
    \begin{equation}\label{eq:sigmadef}
        \sigma_{\gamma,a}^+(s)\coloneqq\max\{\mu\in\Rds:H_\gamma(s,\mu)=a\},\qquad\sigma_{\gamma,a}^-(s)\coloneqq\min\{\mu\in\Rds:H_\gamma(s,\mu)=a\},
    \end{equation}
    with the usual convention that $\sigma_{\gamma,a}^+(s)=-\infty$ and $\sigma_{\gamma,a}^-(s)=\infty$ when $\{\mu\in\Rds:H_\gamma(s,\mu)=a\}$ is empty. It follows from~\cref{eq:condcompH} that
    \begin{equation}\label{eq:sigcomp}
        \sigma_{\widetilde\gamma,a}^+(s)=-\sigma_{\gamma,a}^-(1-s).
    \end{equation}
    Notice that $\{\mu\in\Rds:H_\gamma(s,\mu)=a\}$ is not empty if and only if $a\ge\min\limits_{\mu\in\Rds}H_\gamma(s,\mu)$, thus $\sigma_{\gamma,a}^+(s)\ne-\infty$ for any $s\in[0,1]$ if and only if $a\ge a_\gamma$, where
    \begin{equation}\label{eq:agdef}
        a_\gamma\coloneqq\max_{s\in[0,1]}\min_{\mu\in\Rds}H_\gamma(s,\mu).
    \end{equation}

    \begin{prop}\label{sigcont}\leavevmode
        \begin{myenum}
            \item \label{en:sigcont1} For each $\gamma\in\Ecal$ and $s\in[0,1]$ the function $a\mapsto\sigma_{\gamma,a}^+(s)$ is continuous and increasing in $\mleft[\min\limits_{\mu\in\Rds}H_\gamma(s,\mu),\infty\mright)$.
            \item \label{en:sigcont2} Fixed $\gamma\in\Ecal$, let $[s_1,s_2]\subseteq[0,1]$ and $a\in\Rds$ be such that
            \begin{equation}\label{eq:sigcont.1}
                a\ge\max_{s\in[s_1,s_2]}\min_{\mu\in\Rds}H_\gamma(s,\mu).
            \end{equation}
            Then the function $s\mapsto\sigma_{\gamma,a}^+(s)$ is continuous in $[s_1,s_2]$.
        \end{myenum}
    \end{prop}
    \begin{proof}
        \Cref{en:sigcont1} follows from~\cref{eq:intsublevel}, thus we focus on \cref{en:sigcont2}. We start observing that~\cref{eq:sigcont.1} is equivalent to
        \begin{equation*}
            \sigma_{\gamma,a}^+(s)>-\infty,\qquad \text{for all }s\in[s_1,s_2].
        \end{equation*}
        Let $\{s_n\}_{n\in\Nds}$ be a sequence in $[s_1,s_2]$ converging to a point $\ovs$. The sequence $\mleft\{\sigma_{\gamma,a}^+(s_n)\mright\}_{n\in\Nds}$ is equibounded by~\ref{condcoerc} and therefore converges, up to subsequences, to a $\ovmu\in\Rds$ with $H_\gamma(\ovs,\ovmu)=a$. We have that either $\ovmu=\sigma_{\gamma,a}^+(\ovs)$ or $\ovmu=\sigma_{\gamma,a}^-(\ovs)$. The sequence $\{s_n\}$ is arbitrarily chosen, thus showing that $\ovmu=\sigma_{\gamma,a}^+(\ovs)$ is enough to conclude our proof. If $\sigma_{\gamma,a}^+(\ovs)=\sigma_{\gamma,a}^-(\ovs)$ there is nothing to prove, otherwise \cref{eq:intsublevel} yields
        \begin{equation}\label{eq:sigcont1}
            a>\min_{\mu\in\Rds}H_\gamma(\ovs,\mu).
        \end{equation}
        Thanks to the continuity of the Hamiltonians, \cref{eq:sigcont1} implies the existence of a $\mu_0\in\Rds$ such that
        \begin{equation*}
            \sigma_{\gamma,a}^+(\ovs)>\mu_0>\sigma_{\gamma,a}^-(\ovs)
        \end{equation*}
        and, for any $n$ big enough, $H_\gamma(s_n,\mu_0)<a$, i.e., $\sigma_{\gamma,a}^+(s_n)>\mu_0$. This shows that $\mu_0\le\ovmu$, i.e., $\ovmu=\sigma_{\gamma,a}^+(\ovs)$.
    \end{proof}

    \section{Eikonal Hamilton--Jacobi Equations on Networks}\label{eiksec}

    Here we are interested in equations of the form
    \begin{equation}\label{eq:globeik}\tag{\ensuremath{\Hcal}J\ensuremath{a}}
        \Hcal(x,Du)=a,\qquad \text{on }\Gamma,
    \end{equation}
    thoroughly analyzed in~\cite{SiconolfiSorrentino18}, where $a\in\Rds$. This notation synthetically indicates the family of Hamilton--Jacobi equations
    \begin{equation}\label{eq:eikg}\tag{HJ\textsubscript{\ensuremath{\gamma}}\ensuremath{a}}
        H_\gamma(s,\partial U)=a,\qquad \text{on }[0,1],
    \end{equation}
    for $\gamma$ varying in $\Ecal$.

    Here (sub/super)solutions to the local problem~\cref{eq:eikg} are intended in the viscosity sense, see for instance~\cite{BardiCapuzzo-Dolcetta97} for a comprehensive treatment of viscosity solutions theory. We just recall that, given an open set $\Ocal$ and a continuous function $u:\ovOcal\to\Rds$, a \emph{supertangent} (resp.\ \emph{subtangent}) to $u$ at $x\in\Ocal$ is a viscosity test function from above (resp.\ below). We say that a subtangent $\varphi$ to $u$ at $x\in\partial\Ocal$ is \emph{constrained to $\ovOcal$} if $x$ is a minimizer of $u-\varphi$ in a neighborhood of $x$ intersected with $\ovOcal$.

    \begin{defin}\label{defsol}
        A continuous function $w:\Gamma\to\Rds$ is called a \emph{viscosity subsolution} to~\cref{eq:globeik} if
        \begin{myenum}
            \item $s\mapsto w(\gamma(s))$ is a viscosity subsolution to~\cref{eq:eikg} in $(0,1)$ for any $\gamma\in\Ecal$.
        \end{myenum}
        We say that a continuous function $v:\Gamma\to\Rds$ is a \emph{viscosity supersolution} to~\cref{eq:globeik} if
        \begin{myenum}[resume]
            \item $s\mapsto v(\gamma(s))$ is a viscosity supersolution of~\cref{eq:eikg} in $(0,1)$ for any $\gamma\in\Ecal$;
            \item \label{stateconst} for every vertex $x$ there is at least one arc $\gamma\in\Gamma_x$ such that
            \begin{equation*}
                H_\gamma(1,\partial\varphi(1))\ge a
            \end{equation*}
            for any constrained $C^1$ subtangent $\varphi$ to $v\circ\gamma$ at 1.
        \end{myenum}
        If $u:\Gamma\to\Rds$ is both a viscosity subsolution and supersolution to~\cref{eq:globeik}, we say that $u$ is a \emph{viscosity solution}.
    \end{defin}

    In order to provide a representation formula for solutions to~\cref{eq:globeik}, we extend the support functions defined in~\cref{eq:sigmadef} to the tangent bundle $T\Gamma$ in the following sense: we set, for any $a\in\Rds$, the map $\sigma_a:T\Gamma\to\Rds$ such that
    \begin{mylist}
        \item if $x=\gamma(s)$ for some $\gamma\in\Ecal$ and $s\in(0,1)$ then
        \begin{equation*}
            \sigma_a(x,q)\coloneqq\max\mleft\{\mu\frac{q\dot\gamma(s)}{|\dot\gamma(s)|^2}:\mu\in\Rds,\,H_\gamma(s,\mu)=a\mright\}.
        \end{equation*}
        It is clear that when $\{\mu\in\Rds:H_\gamma(s,\mu)=a\}\ne\emptyset$
        \begin{equation*}
            \sigma_a(x,q)=\max\mleft\{\sigma_{\gamma,a}^+(s)\frac{q\dot\gamma(s)}{|\dot\gamma(s)|^2},\sigma_{\gamma,a}^-(s)\frac{q\dot\gamma(s)}{|\dot\gamma(s)|^2}\mright\},
        \end{equation*}
        otherwise we assume that $\sigma_a(x,q)=-\infty$;
        \item if $x\in\Vbf$ and $q\ne0$ then
        \begin{equation*}
            \sigma_a(x,q)\coloneqq\min\max\mleft\{\mu\frac{q\dot\gamma(1)}{|\dot\gamma(1)|^2}:\mu\in\Rds,\,H_\gamma(1,\mu)=a\mright\},
        \end{equation*}
        where the minimum is taken over the $\gamma\in\Gamma_x$ with $\dot\gamma(1)$ parallel to $q$. We assume that $\sigma_a(x,q)=-\infty$ whenever $\{\mu\in\Rds,\,H_\gamma(1,\mu)=a\}=\emptyset$ for some $\gamma\in\Gamma_x$ with $\dot\gamma(1)$ parallel to $q$;
        \item if $x\in\Vbf$ and $q=0$ then
        \begin{equation*}
            \sigma_a(x,q)\coloneqq0.
        \end{equation*}
    \end{mylist}
    We point out that the case $x\in\Vbf$, $q\ne0$ is more involved because there is a problem to take into account, namely different arcs ending at $x$ could have parallel tangent vectors. In this case we should have a $q\in T_x\Gamma$ parallel to both $\dot\gamma_1(1)$ and $\dot\gamma_2(1)$, where $\gamma_1$ and $\gamma_2$ are different arcs in $\Gamma_x$, i.e.,
    \begin{equation*}
        q=\lambda_1\dot\gamma_1(1)=\lambda_2\dot\gamma_2(1),\qquad \text{for some arcs $\gamma_1\ne\gamma_2$ and scalars }\lambda_1,\lambda_2.
    \end{equation*}
    Notice that, thanks to~\cref{eq:sigcomp}, $\sigma_a$ is a well-defined function in $T\Gamma$.

    Next we set
    \begin{equation}\label{eq:a0def}
        a_0\coloneqq\max_{\gamma\in\Ecal}a_\gamma
    \end{equation}
    and define the \emph{critical value}, or \emph{Mañé critical value}, as
    \begin{equation}\label{eq:critvaldef}
        c\coloneqq\inf\{a\ge a_0: \text{\cref{eq:globeik} admits subsolutions}\}.
    \end{equation}
    The critical value is finite. Indeed, because of \cref{condcoerc}, there is an $a\ge a_0$ large enough so that
    \begin{equation*}
        H_\gamma(s,0)\le a,\qquad \text{for any }s\in[0,1],\,\gamma\in\Ecal,
    \end{equation*}
    i.e., each constant function is a subsolution to~\cref{eq:globeik}. It is also apparent that~\cref{eq:globeik} admits subsolutions whenever $a>c$.

    \begin{rem}\label{critvalchar}
        If $c>a_0$, it is the unique value such that~\hrefc{} (namely the equation~\cref{eq:globeik} with $a=c$) admits solutions in the sense of \cref{defsol}. This is proved in~\cite{SiconolfiSorrentino18} using the following characterization: the critical value is the only $c>a_0$ which satisfies
        \begin{equation}\label{eq:critvalchar}
            \int_0^T\sigma_c\mleft(\xi,\dot\xi\mright)d\tau\ge0,\qquad \text{for all the closed curves }\xi:[0,T]\to\Gamma,
        \end{equation}
        and, for at least one nonconstant closed curve, the inequality above is an identity.\\
        If $c=a_0$ \cref{eq:critvalchar} still holds true, however it is not guaranteed that~\cref{eq:critvalchar} is an identity for some nonconstant closed curve. In~\cite{SiconolfiSorrentino18} this is overcome by requiring the following condition:
        \begin{enumerate}[label=\textbf{(D)},format=\normalfont]
            \item \label{condmina} for any $\gamma\in\Ecal$ with $a_\gamma=c=a_0$ the map $s\mapsto\min\limits_{p\in\Rds}H_\gamma(s,p)$ is constant in $[0,1]$.
        \end{enumerate}
        In some cases, e.g., when we do not know beforehand the value of $c$, condition~\ref{condmina} is too restrictive, thus we do not assume it.\\
        We point out that the results in~\cite{SiconolfiSorrentino18} for~\cref{eq:globeik} when $a>a_0$ still hold true in our setting, since are not affected by condition~\ref{condmina}. In particular, \cref{eq:globeik} does not admit solutions if $a>c$.
    \end{rem}

    Hereafter $c$ will always denote the critical value of the eikonal equation. In addition, (sub/super)solutions to~\hrefc{} will often be referred to as \emph{critical (sub/super)solutions}.

    \begin{rem}
        In~\cite{SiconolfiSorrentino18} the authors give special consideration to \emph{loops}. Namely, a loop is an arc $\gamma$ with $\gamma(0)=\gamma(1)$. They define the set $\Ecal^*$, made up of the loops in $\Ecal$, and $c_\gamma$ as the minimum $a\ge a_\gamma$ such that~\cref{eq:eikg} admits periodic subsolutions. We stress out that we do not assume any periodicity on $H_\gamma$ when $\gamma$ is a loop. They further set
        \begin{equation*}
            a_0\coloneqq\max\mleft\{\max_{\gamma\in\Ecal\setminus\Ecal^*}a_\gamma,\max_{\gamma\in\Ecal^*}c_\gamma\mright\},
        \end{equation*}
        which is clearly bigger than or equal to the $a_0$ defined in~\cref{eq:a0def}. However, the definition of critical value~\cref{eq:critvaldef} remain the same regardless of the $a_0$ used: it is indeed apparent that if \cref{eq:globeik} has a subsolution $w$ and $\gamma$ is a loop then $w\circ\gamma$ is a periodic subsolution to \cref{eq:eikg}, i.e., $a\ge c_\gamma$. In view of this fact and since loops do not require special attention in our analysis, throughout this paper we will employ the definition of $a_0$ given in~\cref{eq:a0def}.
    \end{rem}

    We conclude this \lcnamecref{eiksec} by defining the semidistance on $\Gamma$
    \begin{equation}\label{eq:semidist}
        S_a(y,x)\coloneqq\inf\mleft\{\int_0^T\sigma_a\mleft(\xi,\dot\xi\mright)d\tau: \text{$\xi:[0,T]\to\Gamma$ is a curve from $y$ to }x\mright\},
    \end{equation}
    whose importance is highlighted by the next \namecref{subsolchar}.

    \begin{prop}\label{subsolchar}
        A continuous function $w:\Gamma\to\Rds$ is a subsolution to~\cref{eq:globeik} if and only if
        \begin{equation*}
            w(x)-w(y)\le S_a(y,x),\qquad \text{for any }x,y\in\Gamma.
        \end{equation*}
    \end{prop}
    \begin{proof}
        See~\cite{SiconolfiSorrentino18}.
    \end{proof}

    The proofs of the next results, which characterize the semidistance~\cref{eq:semidist}, are given in \cref{optcurve}.

    \begin{defin}\label{subarcdef}
        Let $\zeta:[0,T]\to\Gamma$ be a simple curve. We say that $\zeta$ is a \emph{sub-arc} if there exist an arc $\gamma$ and a curve $\eta:[0,T]\to[0,1]$ with $\dot\eta=1$ a.e.\ such that $\zeta\equiv\gamma\circ\eta$.
    \end{defin}

    \begin{restatable}{prop}{mincurve}\label{mincurve}
        Let $x,y\in\Gamma$ and $a\ge c$. The infimum in~\cref{eq:semidist} is a minimum and,
        \begin{myenum}
            \item \label{en:mincurve1} $S_a(x,x)=0$;
            \item \label{en:mincurve2} if $x\ne y$, there is a simple curve $\zeta:[0,T]\to\Gamma$, which is a concatenation of sub-arcs linking $y$ and $x$, such that
            \begin{equation*}
                S_a(y,x)=\int_0^T\sigma_a\mleft(\zeta,\dot\zeta\mright)d\tau.
            \end{equation*}
        \end{myenum}
    \end{restatable}

    \begin{restatable}{prop}{sdiscont}\label{sdiscont}
        For any $a\ge c$ the map
        \begin{equation*}
            \Gamma^2\ni(y,x)\longrightarrow S_a(y,x)
        \end{equation*}
        is Lipschitz continuous with respect to the geodesic distance $d_\Gamma$.
    \end{restatable}

    \section{Aubry Set and Representation Formula}\label{aubsec}

    The following set, whose definition is deeply related to the critical value $c$, is crucial for our analysis:

    \begin{defin}\label{defaubry}
        We call \emph{Aubry set} on $\Gamma$, the closed set $\Acal_\Gamma$ made up of
        \begin{myenum}
            \item \label{en:defaubry1} the $x\in\Gamma$ incident to a closed simple curve $\xi:[0,T]\to\Gamma$ with $\int_0^T\sigma_c\mleft(\xi,\dot\xi\mright)d\tau=0$;
            \item \label{en:defaubry2} the $x=\gamma(s)$ with $\gamma\in\Ecal$ and $s\in[0,1]$ such that $\sigma_{\gamma,c}^+(s)=\sigma_{\gamma,c}^-(s)$.
        \end{myenum}
    \end{defin}

    We point out that, if $c=a_0$, then there must be a $\gamma\in\Ecal$ with $a_\gamma=c$ and consequently, by \cref{condqconv}, there is at least one point $x$ as in \cref{defaubry}\ref{en:defaubry2}. If instead $c>a_0$, we know from~\cref{eq:critvalchar} that there exists at least one closed curve as in \cref{defaubry}\ref{en:defaubry1}. This shows that the Aubry set is always nonempty.

    \begin{rem}\label{diffAub}
        A different notion of Aubry set has been given in~\cite{SiconolfiSorrentino18}, which is solely made up of the $x\in\Gamma$ incident to a closed curve $\xi:[0,T]\to\Gamma$ with $\int_0^T\sigma_c\mleft(\xi,\dot\xi\mright)d\tau=0$ and a.e.\ non-vanishing derivative. Under their assumptions, such set consists of the support of a collection of arcs. Conversely, in our case, the Aubry set $\Acal_\Gamma$ could even be made up of a countable number of disjointed pieces of arcs or just a single point.
    \end{rem}

    Let us assume that there exist a $\gamma\in\Ecal$ and a closed interval $[s_1,s_2]\subseteq[0,1]$ such that
    \begin{equation}\label{eq:flatlevel}
        \sigma_{\gamma,c}^+(s)=\sigma_{\gamma,c}^-(s),\qquad \text{for any }s\in[s_1,s_2].
    \end{equation}
    We define the curves
    \begin{equation*}
        \eta(t)\coloneqq\mleft\{
        \begin{aligned}
            &t+s_1,&& \text{if }t\in[0,s_2-s_1),\\
            &2s_2-s_1-t,&& \text{if }t\in[s_2-s_1,2(s_2-s_1)],
        \end{aligned}
        \mright.
    \end{equation*}
    and $\xi\coloneqq\gamma\circ\eta$, then
    \begin{equation*}
        \int_0^{2(s_2-s_1)}\sigma_c\mleft(\xi,\dot\xi\mright)d\tau=\int_0^{s_2-s_1}\sigma_{\gamma,c}^+(r)dr-\int_{s_2-s_1}^{2(s_2-s_1)}\sigma_{\gamma,c}^-(r)dr=0.
    \end{equation*}
    It is therefore apparent that a point in the Aubry set is either incident to a closed curve $\xi:[0,T]\to\Gamma$ with $\int_0^T\sigma_c\mleft(\xi,\dot\xi\mright)d\tau=0$ and a.e.\ non-vanishing derivative, or an isolated point of $\Acal_\Gamma$. Moreover, a maximal interval $I$ such that $\sigma_{\gamma,c}^+=\sigma_{\gamma,c}^-$ on $I$ must be closed by \cref{sigcont}. These facts suggest the \namecref{statclass} below.

    \begin{defin}\label{statclass}
        The Aubry set is partitioned into \emph{static classes}, defined as the singletons containing the isolated points of $\Acal_\Gamma$ or the equivalence classes with respect to the relation
        \begin{equation*}
            \mleft\{
            \begin{gathered}
                \text{$x,y\in\Gamma$: $x$ and $y$ are incident to a closed curve $\xi:[0,T]\to\Gamma$ with }\int_0^T\sigma_c\mleft(\xi,\dot\xi\mright)d\tau=0\\
                \text{and a.e.\ non-vanishing derivative}
            \end{gathered}
            \mright\}.
        \end{equation*}
    \end{defin}

    \begin{rem}
        It easily follows from~\cref{eq:agdef,eq:sigmadef,eq:a0def,eq:critvaldef}, that \cref{en:defaubry2} in \cref{defaubry} is relevant for the definition of $\Acal_\Gamma$ only when $c=a_0$. Moreover, if we assume condition~\ref{condmina}, \cref{eq:flatlevel} holds true, with $s_1=0$ and $s_2=1$, for any $\gamma$ with $a_\gamma=a_0=c$. These facts yield that if~\ref{condmina} holds or $c>a_0$, then every point in the Aubry set is incident to a closed curve $\xi:[0,T]\to\Gamma$ with $\int_0^T\sigma_c\mleft(\xi,\dot\xi\mright)d\tau=0$ and a.e.\ non-vanishing derivative. Consequently, see \cref{diffAub}, $\Acal_\Gamma$ corresponds to the Aubry set defined in~\cite{SiconolfiSorrentino18}.
    \end{rem}

    We exploit the properties of the Aubry set to further characterize critical subsolutions.

    \begin{lem}\label{subsolstatclass}
        If $\Ccal$ is a static class of $\Acal_\Gamma$ and $w$ is a subsolution to~\hrefc{} then
        \begin{equation*}
            w(x)=w(y)+S_c(y,x),\qquad \text{for any }x,y\in\Ccal.
        \end{equation*}
        Furthermore, if $\zeta:[0,T]\to\Gamma$ is a curve with $\int_0^T\sigma_c\mleft(\zeta,\dot\zeta\mright)d\tau=0$, then
        \begin{equation}\label{eq:subsolstatclass.1}
            S_c\mleft(\zeta(t),\zeta\mleft(t'\mright)\mright)=\int_t^{t'}\sigma_c\mleft(\zeta,\dot\zeta\mright)d\tau
        \end{equation}
        for any $0\le t\le t'\le T$.
    \end{lem}
    \begin{proof}
        If $\Ccal=\{x\}$ our claim is trivially true, see \cref{mincurve}. We thus assume that this is not the case. We argue by contradiction, assuming that there exist $x$, $y$ in a static class $\Ccal$ such that
        \begin{equation}\label{eq:subsolstatclass1}
            w(x)\ne w(y)+S_c(y,x).
        \end{equation}
        Thanks to \cref{statclass}, we have two curves $\zeta_1:[0,T_1]\to\Ccal$ and $\zeta_2:[0,T_2]\to\Ccal$ with $\zeta_1(0)=\zeta_2(T_2)=x$, $\zeta_1(T_1)=\zeta_2(0)=y$ and a.e.\ non-vanishing derivative such that, setting $\zeta\coloneqq\zeta_1*\zeta_2$ and $T\coloneqq T_1+T_2$,
        \begin{equation}\label{eq:subsolstatclass2}
            \int_0^{T_1}\sigma_c\mleft(\zeta_1,\dot\zeta_1\mright)d\tau+\int_0^{T_2}\sigma_c\mleft(\zeta_2,\dot\zeta_2\mright)d\tau=\int_0^T\sigma_c\mleft(\zeta,\dot\zeta\mright)d\tau=0.
        \end{equation}
        It follows from \cref{subsolchar} and~\cref{eq:subsolstatclass1} that
        \begin{equation*}
            w(x)-w(y)<S_c(y,x)\le\int_0^{T_1}\sigma_c\mleft(\zeta_1,\dot\zeta_1\mright)d\tau,\qquad w(y)-w(x)\le S_c(x,y)\le\int_0^{T_2}\sigma_c\mleft(\zeta_2,\dot\zeta_2\mright)d\tau,
        \end{equation*}
        thus
        \begin{equation*}
            0=w(x)-w(y)+w(y)-w(x)<S_c(y,x)+S_c(x,y)\le\int_0^T\sigma_c\mleft(\zeta,\dot\zeta\mright)d\tau,
        \end{equation*}
        in contradiction with~\cref{eq:subsolstatclass2}. The same arguments also prove~\cref{eq:subsolstatclass.1}.
    \end{proof}

    \begin{prop}\label{subsolonaub}
        Any subsolution to~\hrefc{} is also a solution in $\Acal_\Gamma$.
    \end{prop}
    \begin{proof}
        We start showing that, given a subsolution $w$, $\gamma\in\Ecal$ and $s\in(0,1)$ such that $\gamma(s)\in\Acal_\Gamma\setminus\Vbf$, $w\circ\gamma$ is differentiable at $s$ and
        \begin{equation}\label{eq:subsolonaub1}
            H_\gamma(s,D(w\circ\gamma)(s))=c.
        \end{equation}
        It is well known, see for instance~\cite{BardiCapuzzo-Dolcetta97}, that $D(w\circ\gamma)(s)=D\varphi(s)$ for any $C^1$ subtangent $\varphi$ to $w\circ\gamma$ at $s$, \cref{eq:subsolonaub1} thus yields that $w\circ\gamma$ is a solution to~\ghrefc{} at $s$ and, consequently, $w$ is a solution to~\hrefc{} at $\gamma(s)$. Since $w$, $\gamma$ and $s$ are arbitrary, this will prove that any critical subsolution is a solution to~\hrefc{} in $\Acal_\Gamma\setminus\Vbf$.\\
        We first assume that $\gamma$ and $s$ are as in \cref{defaubry}\ref{en:defaubry2}, namely
        \begin{equation}\label{eq:subsolonaub2}
            \sigma_{\gamma,c}^+(s)=\sigma_{\gamma,c}^-(s).
        \end{equation}
        Thanks to \cref{subsolchar} we have, for $h>0$ small enough,
        \begin{align*}
            w(\gamma(s+h))-w(\gamma(s))&\le S_c(\gamma(s),\gamma(s+h))\le\int_s^{s+h}\sigma^+_{\gamma,c}(r)dr,\\
            w(\gamma(s+h))-w(\gamma(s))&\ge-S_c(\gamma(s+h),\gamma(s))\ge-\int_s^{s+h}\sigma^+_{\wtgamma,c}(1-r)dr=\int^{s+h}_s\sigma^-_{\gamma,c}(r)dr,
        \end{align*}
        which implies
        \begin{equation*}
            \sigma^-_{\gamma,c}(s)=\lim_{h\to0^+}\frac{w(\gamma(s+h))-w(\gamma(s))}h=\sigma^+_{\gamma,c}(s).
        \end{equation*}
        Similarly, we can show that the left derivative of $w\circ\gamma$ at $s$ is $\sigma^+_{\gamma,c}(s)$, therefore $D(w\circ\gamma)(s)=\sigma^+_{\gamma,c}(s)$ and~\cref{eq:subsolonaub1} holds true.\\
        Next we assume that~\cref{eq:subsolonaub2} is false, then \cref{defaubry,lowcost,closedcurve} yield that, possibly replacing $\gamma$ with $\wtgamma$, there is a closed simple curve $\zeta:[0,T]\to\Gamma$, which is a concatenation of sub-arcs, with $\zeta(0)=\gamma(s)$ and $\int_0^T\sigma_c\mleft(\zeta,\dot\zeta\mright)d\tau=0$. We get from \cref{subsolstatclass} that, for any $h>0$ small enough,
        \begin{equation*}
            w(\gamma(s+h))-w(\gamma(s))=S_c(\gamma(s),\gamma(s+h))=\int_0^h\sigma_c\mleft(\zeta,\dot\zeta\mright)d\tau=\int_s^{s+h}\sigma^+_{\gamma,c}(r)dr
        \end{equation*}
        and
        \begin{equation*}
            w(\gamma(s))-w(\gamma(s-h))=S_c(\gamma(s-h),\gamma(s))=\int_{T-h}^T\sigma_c\mleft(\zeta,\dot\zeta\mright)d\tau=\int^s_{s-h}\sigma^+_{\gamma,c}(r)dr.
        \end{equation*}
        It follows that $w\circ\gamma$ has both left and right derivative at $s$ equal to $\sigma_{\gamma,c}^+(s)$, i.e., it is differentiable at $s$ with $D(w\circ\gamma)(s)=\sigma_{\gamma,c}^+(s)$, thus~\cref{eq:subsolonaub1} holds true.\\
        Now let $x$ be a vertex in $\Acal_\Gamma$. Arguing as in the previous part of the proof we get that there is an arc $\gamma$ such that $x=\gamma(1)$ and $w\circ\gamma$ is left differentiable at 1 with left derivative equal to $\sigma_{\gamma,c}^+(1)$. Let $\varphi$ be a constrained $C^1$ subtangent to $w\circ\gamma$ at 1, then
        \begin{equation*}
            \partial\varphi(1)=\lim_{h\to0^+}\frac{\varphi(1)-\varphi(1-h)}h\ge\lim_{h\to0^+}\frac{w(\gamma(1))-w(\gamma(1-h))}h=\sigma^+_{\gamma,c}(1),
        \end{equation*}
        which proves that $w$ satisfies \labelcref{stateconst} in \cref{defsol}. This concludes our proof since $w$ and $x$ are arbitrary.
    \end{proof}

    \begin{rem}
        The proof of \cref{subsolonaub} shows a regularity property of critical subsolutions: if $x\in\Acal_\Gamma\cap\Vbf$, even if it is an isolated point of the Aubry set, there exist $\gamma\in\Ecal$ and $s\in(0,1)$ such that $\gamma(s)=x$ and
        \begin{equation*}
            D(w\circ\gamma)(s)=\sigma_{\gamma,c}^+(s).
        \end{equation*}
        Moreover, \cref{subsolstatclass} proves that critical subsolutions are uniquely determined in a static class by its value at a single point, i.e., they differ by a constant. It follows that all the subsolutions to~\hrefc{} possess the same differential in $\Acal_\Gamma\setminus\Vbf$ (see \cref{diffdef}), extending \cite[Theorem 7.5]{SiconolfiSorrentino18} to our case.
    \end{rem}

    The next \namecref{maxsubsol}, which extends the existence result given in~\cite{SiconolfiSorrentino18} to our setting, is the main connection between the Aubry set and critical (sub)solutions.

    \begin{theo}\label{maxsubsol}
        Let $\Gamma'$ be a closed subset of $\Gamma$, $g:\Gamma'\to\Rds$ be a continuous function and define
        \begin{equation}\label{eq:maxsubsol.1}
            u(x)\coloneqq\min_{y\in\Gamma'}(g(y)+S_c(y,x)),\qquad \text{for }x\in\Gamma.
        \end{equation}
        Then $u$ is the maximal subsolution to~\hrefc{} not exceeding $g$ on $\Gamma'$ and a solution in $\Gamma\setminus(\Gamma'\setminus\Acal_\Gamma)$.
    \end{theo}
    \begin{proof}
        First observe that \cref{sdiscont} justifies the minimum in~\cref{eq:maxsubsol.1} and yields the continuity of $u$. In addition, $u$ is a subsolution because of \cref{subsolchar}. The maximality is proved arguing by contradiction: let $x\in\Gamma$ and $w$ be a critical subsolution not exceeding $g$ on $\Gamma'$ such that $u(x)<w(x)$. For any $y$ optimal to $u(x)$ in~\cref{eq:maxsubsol.1} we have
        \begin{equation*}
            w(y)+S_c(y,x)\le g(y)+S_c(y,x)=u(x)<w(x),
        \end{equation*}
        i.e., $S_c(y,x)<w(x)-w(y)$, in contradiction with \cref{subsolchar}.\\
        It only remains to prove that $u$ is a critical supersolution at $x$ whenever $x\in\Gamma\setminus(\Gamma'\setminus\Acal_\Gamma)$. We distinguish three cases:
        \begin{myenum}
            \item \label{en:maxsubsol1} $x\in\Acal_\Gamma$;
            \item \label{en:maxsubsol2} $x\in\Gamma\setminus(\Gamma'\cup\Vbf)$;
            \item \label{en:maxsubsol3} $x\in\Vbf\setminus\Gamma'$.
        \end{myenum}
        We will focus on case~\ref{en:maxsubsol3}. The same arguments will also prove case~\ref{en:maxsubsol2}, while case~\ref{en:maxsubsol1} has already been proved in \cref{subsolonaub}.\\
        By~\cref{eq:maxsubsol.1} and \cref{mincurve}, there is an $y\in\Gamma'$ and a simple curve $\zeta:[0,T]\to\Gamma$, which is a concatenation of sub-arcs, such that
        \begin{equation*}
            u(x)=g(y)+S_c(y,x)=g(y)+\int_0^T\sigma_c\mleft(\zeta,\dot\zeta\mright)d\tau.
        \end{equation*}
        The optimality of $\zeta$ yields that there is a $\gamma\in\Gamma_x$ such that, for any $h$ small enough,
        \begin{align*}
            u(x)&=u(\gamma(1))=u(\gamma(1-h))+S_c(\gamma(1-h),\gamma(1))=u(\gamma(1-h))+\int_{T-h}^T\sigma_c\mleft(\zeta,\dot\zeta\mright)d\tau\\
            &=u(\gamma(1-h))+\int_{1-h}^1\sigma_{\gamma,c}^+(r)dr.
        \end{align*}
        Let $\varphi$ be any constrained $C^1$ subtangent to $u\circ\gamma$ at 1, then
        \begin{equation*}
            \partial\varphi(1)=\lim_{h\to0^+}\frac{\varphi(1)-\varphi(1-h)}h\ge\lim_{h\to0^+}\frac{u(\gamma(1))-u(\gamma(1-h))}h=\sigma^+_{\gamma,c}(1),
        \end{equation*}
        which proves that $u$ satisfies \labelcref{stateconst} in \cref{defsol}, i.e., it is a critical supersolution at $x$.
    \end{proof}

    The next \namecref{critvalunique} is a simple consequence of \cref{critvalchar,maxsubsol}:

    \begin{cor}\label{critvalunique}
        If \cref{eq:globeik} admits solutions, then $a=c$.
    \end{cor}

    Taking into account \cref{subsolchar,maxsubsol}, we call \emph{admissible trace} any continuous real function $g$ defined on a subset $\Gamma'\subseteq\Gamma$ such that
    \begin{equation}\label{eq:trace}
        g(x)-g(y)\le S_c(y,x),\qquad \text{for any }x,y\in\Gamma'.
    \end{equation}
    We point out that every subsolution is also an admissible trace. The next \namecref{hopflax} is a trivial consequence of \cref{maxsubsol}.

    \begin{cor}\label{hopflax}
        Let $\Gamma'$ be a closed subset of $\Gamma$. Given an admissible trace $g:\Gamma'\to\Rds$, in the sense of~\cref{eq:trace}, the function
        \begin{equation}\label{eq:hopflax.1}
            u(x)\coloneqq\min_{y\in\Gamma'}(g(y)+S_c(y,x)),\qquad \text{for }x\in\Gamma,
        \end{equation}
        is the maximal subsolution to~\hrefc{} agreeing with $g$ on $\Gamma'$ and a solution in $\Gamma\setminus(\Gamma'\setminus\Acal_\Gamma)$.
    \end{cor}

    We point out that in~\cite{SiconolfiSorrentino18} it is also proved that, when $\Gamma'=\Acal_\Gamma$, \cref{eq:hopflax.1} is the \emph{unique} critical solution agreeing with $g$ on $\Acal_\Gamma$. Their proof heavily relies on condition~\ref{condmina}, thus it does not apply to our case. We will show later that, under our assumptions, such solution is still unique.

    \section{Comparison Result for the Eikonal Equation}\label{compsec}

    This \lcnamecref{compsec} is devoted to prove a comparison principle for the eikonal problem~\hrefc{} where the Aubry set plays, in a sense, the role of a hidden boundary. To our knowledge there is no previous comparison result between super and subsolutions for such equation.

    Some preliminary results are needed. The first one is self-evident.

    \begin{lem}\label{strictsubsol}
        Let $w_\gamma$ be a subsolution to~\cref{eq:eikg} in $(s_1,s_2)\subseteq(0,1)$ and assume that
        \begin{equation*}
            \sigma_{\gamma,a}^+(s)>\sigma_{\gamma,a}^-(s),\qquad \text{for all }s\in(s_1,s_2).
        \end{equation*}
        Fixed $0<\delta<s_2-s_1$, there is a sequence $\{w_n\}_{n\in\Nds}\subset C^1([s_1+\delta,s_2-\delta])$ of subsolutions uniformly converging to $w_\gamma$ on $[s_1+\delta,s_2-\delta]$ such that
        \begin{equation*}
            H_\gamma(s,Dw_n(s))<a,\qquad \text{for every }s\in[s_1+\delta,s_2-\delta],\,n\in\Nds.
        \end{equation*}
    \end{lem}

    \begin{lem}\label{locmaxsubsol}
        Let $U$ be a solution to~\cref{eq:eikg} in $(\ovs,1)\subseteq(0,1)$, continuously extended up to $[\ovs,1]$. If
        \begin{equation*}
            \sigma_{\gamma,a}^+(s)>\sigma_{\gamma,a}^-(s),\qquad \text{for all }s\in(\ovs,1],
        \end{equation*}
        and
        \begin{equation*}
            H_\gamma(1,\partial\varphi(1))\ge a,\qquad \text{for any constrained $C^1$ supertangent $\varphi$ to $U$ at 1},
        \end{equation*}
        then
        \begin{equation*}
            U(s)=U(\ovs)+\int_{\ovs}^s\sigma_{\gamma,a}^+(r)dr,\qquad \text{for each }s\in[\ovs,1].
        \end{equation*}
    \end{lem}
    \begin{proof}
        For any $\delta>0$ small enough, \cite[Proposition 5.6]{SiconolfiSorrentino18} yields that
        \begin{equation*}
            U(s)=U(\ovs+\delta)+\int_{\ovs+\delta}^s\sigma_{\gamma,a}^+(r)dr,\qquad \text{for each }s\in[\ovs+\delta,1].
        \end{equation*}
        Since $U$ is continuous and $\delta$ is arbitrary, this proves our claim.
    \end{proof}

    \begin{theo}[Comparison Principle]\label{eikcomp}
        Let $v$ be a continuous supersolution to~\hrefc{} in $\Gamma\setminus\Gamma'$, where $\Gamma'$ is a closed subset of $\Gamma$ containing $\Acal_\Gamma$, and $w$ be a critical subsolution. If $v\ge w$ on $\Gamma'$, then $v\ge w$ on $\Gamma$.
    \end{theo}

    The natural choice of $\Gamma'$ in the above statement is clearly $\Acal_\Gamma$. There are however situations where considering a general $\Gamma'\supset\Acal_\Gamma$ makes more sense, such as numerical applications in which the Aubry set cannot be explicitly represented (see \cref{diffAub}). See also~\cite{Pozza23}, where the Aubry set is, in a certain sense, extended by the asymptotic behavior of the solutions to a related time-dependent problem.

    \begin{proof}
        First we define
        \begin{equation*}
            u(x)\coloneqq\min_{y\in\Gamma'}(w(y)+S_c(y,x)),\qquad \text{for }x\in\Gamma,
        \end{equation*}
        which by \cref{hopflax} is both a critical solution in $\Gamma\setminus\Gamma'$ and the maximal subsolution in $\Gamma$ agreeing with $w$ on $\Gamma'$. We will show that there is a minimizer $z\in\Gamma'$ to $v-u$, which implies
        \begin{equation*}
            v(z)-w(z)=v(z)-u(z)\le v(x)-u(x)\le v(x)-w(x),\qquad \text{for any }x\in\Gamma,
        \end{equation*}
        i.e., $v-w$ achieves its minimum in $\Gamma'$ proving our claim. We argue by contradiction, assuming that all the minimizer of $v-u$ are outside $\Gamma'$.\\
        Let $\gamma\in\Ecal$ and $(s_1,s_2)\subset[0,1]$ be such that $\gamma((s_1,s_2))\cap\Gamma'=\emptyset$, then we have
        \begin{equation*}
            \sigma_{\gamma,c}^+(s)>\sigma_{\gamma,c}^-(s),\qquad \text{for any }s\in(s_1,s_2).
        \end{equation*}
        Fixed $0<\delta<s_2-s_1$ we get, thanks to \cref{strictsubsol}, a sequence $\{u_{\gamma,n}\}_{n\in\Nds}\subset C^1([s_1+\delta,s_2-\delta])$ uniformly converging to $u\circ\gamma$ on $[s_1+\delta,s_2-\delta]$ as $n\to\infty$ such that
        \begin{equation}\label{eq:eikcomp1}
            H_\gamma(s,\partial u_{\gamma,n}(s))<c,\qquad \text{for every }s\in[s_1+\delta,s_2-\delta],\,n\in\Nds.
        \end{equation}
        If $v\circ\gamma-u_{\gamma,n}$ has a local minimum at $s\in(s_1+\delta,s_2-\delta)$ then $u_{\gamma,n}$ is a $C^1$ subtangent to $v\circ\gamma$ at $s$, which is in contradiction with~\cref{eq:eikcomp1} and the supersolution property of $v\circ\gamma$. Therefore, at least one between $s_1+\delta$ and $s_2-\delta$ is a minimizer of $v\circ\gamma-u_{\gamma,n}$ in $[s_1+\delta,s_2-\delta]$. Since $\delta$ is arbitrary and $u_{\gamma,n}$ uniformly converges to $u\circ\gamma$, this shows that $v\circ\gamma-u\circ\gamma$ restricted to $[s_1,s_2]$ achieves its minimum at $s_1$ or $s_2$. Moreover, also $s_1$, $s_2$ and $\gamma$ are arbitrarily chosen, thus our assumptions yield that there is a minimizer $\ovx$ of $v-u$ contained in $\Vbf\setminus\Gamma'$, which we assume nonempty. We choose a $\gamma_1\in\Gamma_{\ovx}$ so that $v$, $\ovx$ and $\gamma_1$ satisfy \labelcref{stateconst} in \cref{defsol}. Notice that 1 is a minimizer for $v\circ\gamma_1-u\circ\gamma_1$, thereby a constrained subtangent to $u\circ\gamma$ at 1 is also a constrained subtangent to $v\circ\gamma$ at 1, namely $u$, $\ovx$ and $\gamma_1$ satisfy \labelcref{stateconst} in \cref{defsol}. We further know from the definition of Aubry set and \cref{sigcont} that there is an $\ovs\in[0,1)$ such that
        \begin{equation}\label{eq:eikcomp2}
            \sigma_{\gamma_1,c}^+(s)>\sigma_{\gamma_1,c}^-(s),\qquad \text{for all }s\in(\ovs,1].
        \end{equation}
        Hence $u\circ\gamma_1$ satisfies the hypotheses of \cref{locmaxsubsol}, yielding
        \begin{equation*}
            u(\gamma_1(s))=u(\gamma_1(\ovs))+\int_{\ovs}^s\sigma_{\gamma_1,c}^+(\tau)d\tau,\qquad \text{for any }s\in[\ovs,1].
        \end{equation*}
        Exploiting~\cref{eq:eikcomp2,condqconv,sigcont}, we set an $n_0>0$ and an infinitesimal sequence $\{\delta_n\}_{n\ge n_0}$ such that
        \begin{equation*}
            \sigma_{\gamma_1,c-\frac1n}^+(s)>\sigma_{\gamma_1,c-\frac1n}^-(s),\qquad \text{for all }n\ge n_0,\,s\in[\ovs+\delta_n,1].
        \end{equation*}
        We further define, for every $n\ge n_0$,
        \begin{equation*}
            U_n(s)\coloneqq u(\gamma_1(\ovs+\delta_n))+\int_{\ovs+\delta_n}^s\sigma_{\gamma_1,c-\frac1n}^+(\tau)d\tau,\qquad \text{for }s\in[\ovs+\delta_n,1],
        \end{equation*}
        which is a $C^1([\ovs+\delta_n,1])$ function such that
        \begin{equation}\label{eq:eikcomp3}
            H_{\gamma_1}(s,\partial U_n(s))<c,\qquad \text{for any }n\ge n_0,\,s\in[\ovs+\delta_n,1].
        \end{equation}
        As in the previous step, the definition of supersolution implies that the only local minimum of $v\circ\gamma_1-U_n$ is achieved at $\ovs+\delta_n$. Indeed, if $s\in(\ovs+\delta_n,1]$ is a local minimizer, then $U_n$ is a $C^1$ subtangent to $v\circ\gamma_1$ at $s$ (a constrained subtangent if $s=1$) satisfying~\cref{eq:eikcomp3}, in contradiction with our assumptions. It follows that $v\circ\gamma_1-U_n$ is increasing in $[\ovs+\delta_n,1]$ and in particular
        \begin{equation*}
            v(\gamma_1(\ovs+\delta_n))-U_n(\ovs+\delta_n)\le v(\gamma_1(s))-U_n(s)\le v(\gamma_1(1))-U_n(1),\qquad \text{for any }s\in[\ovs+\delta_n,1].
        \end{equation*}
        By definition $U_n$ locally uniformly converges on $(\ovs,1]$ to $u\circ\gamma_1$ as $n$ tends to $\infty$, while 1 is a minimizer of $v\circ\gamma_1-u\circ\gamma_1$. Consequently, the previous inequality yields that $v\circ\gamma_1-u\circ\gamma_1$ is constant on $[\ovs,1]$ and, since the minimizers of $v-u$ are not in $\Gamma'$, $\gamma_1([\ovs,1])\cap\Gamma'=\emptyset$. We point out that $\Gamma'$ is closed, therefore we can assume that $\ovs$ in~\cref{eq:eikcomp2} is equal to 0, i.e.,
        \begin{equation*}
            \gamma_1([0,1])\cap\Gamma'=\emptyset\qquad \text{and}\qquad u(\gamma_1(s))=u(\gamma_1(0))+\int_0^s\sigma_{\gamma_1,c}^+(\tau)d\tau,\quad \text{for any }s\in[0,1].
        \end{equation*}
        Next we choose another arc $\gamma_2\in\Gamma_{\gamma_1(0)}$ such that \labelcref{stateconst} in \cref{defsol} holds for $v$, $\gamma_2$ and $\gamma_1(0)=\gamma_2(1)$. Arguing as before, we get
        \begin{equation*}
            \gamma_2([0,1])\cap\Gamma'=\emptyset\qquad \text{and}\qquad u(\gamma_2(s))=u(\gamma_2(0))+\int_0^s\sigma_{\gamma_2,c}^+(\tau)d\tau,\quad \text{for any }s\in[0,1].
        \end{equation*}
        Iterating this procedure $n$ times, we get a concatenation of arcs $\xi\coloneqq\gamma_n*\dotsb*\gamma_1$ with
        \begin{equation*}
            \xi([0,n])\cap\Gamma'=\emptyset
        \end{equation*}
        and
        \begin{equation*}
            u(\xi(t))=u(\xi(0))+\int_0^t\sigma_c\mleft(\xi,\dot\xi\mright)d\tau,\qquad \text{for any }t\in[0,n].
        \end{equation*}
        The arcs are finite, therefore after a finite number of iterations we get a closed curve $\zeta\coloneqq\gamma'_k*\dotsb*\gamma_1'$ such that
        \begin{equation}\label{eq:eikcomp4}
            \zeta([0,k])\cap\Gamma'=\emptyset
        \end{equation}
        and
        \begin{equation*}
            u(\zeta(t))=u(\zeta(0))+\int_0^t\sigma_c\mleft(\zeta,\dot\zeta\mright)d\tau,\qquad \text{for any }t\in[0,k].
        \end{equation*}
        Finally we have
        \begin{equation*}
            \int_0^k\sigma_c\mleft(\zeta,\dot\zeta\mright)d\tau=u(\zeta(k))-u(\zeta(0))=0,
        \end{equation*}
        which is in contradiction with~\cref{eq:eikcomp4} because $\Gamma'$ contains the Aubry set.
    \end{proof}

    Notice that this comparison result can be carried effortlessly to the supercritical case. Indeed, keeping in mind that when $a>c$
    \begin{equation*}
        \sigma_{\gamma,a}^+(s)>\sigma_{\gamma,a}^-(s),\qquad \text{for all }\gamma\in\Ecal,\,s\in[0,1],
    \end{equation*}
    and, for any nonconstant closed curve $\xi:[0,T]\to\Gamma$, (see \cref{critvalchar})
    \begin{equation*}
        \int_0^T\sigma_a\mleft(\xi,\dot\xi\mright)d\tau>0,
    \end{equation*}
    the proof of \cref{eikcomp}, with straightforward modifications, also proves the next \namecref{supcritcomp}.

    \begin{theo}\label{supcritcomp}
        Let $a>c$ and $\Gamma'$ be a closed subset of $\Gamma$. If $v$ and $w$ are a continuous supersolution in $\Gamma\setminus\Gamma'$ and a subsolution in $\Gamma$ to~\cref{eq:globeik}, respectively, with $v\ge w$ on $\Gamma'$, then $v\ge w$ on $\Gamma$.
    \end{theo}

    Combining \cref{eikcomp,hopflax}, we obtain that the Hopf--Lax type formula~\cref{eq:hopflax.1} provides a representation formula for the unique critical solution to~\hrefc{} assuming a prescribed value on the Aubry set:

    \begin{theo}[Existence and uniqueness of solutions]
        Given an admissible trace $g$, in the sense of~\cref{eq:trace}, the function
        \begin{equation*}
            u(x)\coloneqq\min_{y\in\Acal_\Gamma}(g(y)+S_c(y,x)),\qquad \text{for }x\in\Gamma,
        \end{equation*}
        is the unique solution to~\hrefc{} agreeing with $g$ on $\Acal_\Gamma$.
    \end{theo}

    \appendix

    \section{Optimal Curves}\label{optcurve}

    The optimal curves for the semidistance~\cref{eq:semidist} play a central role in our analysis. This \lcnamecref{optcurve} is therefore devoted to the characterization of such curves.

    We start with a result about curves whose support is contained in an arc of the network.

    \begin{lem}\label{invarccomp}
        \emph{\cite[Lemma 3.2]{PozzaSiconolfi23}} For any given arc $\gamma$ and curve $\xi:[0,T]\to\gamma([0,1])$, the function
        \begin{equation*}
            \gamma^{-1}\circ\xi:[0,T]\to[0,1]
        \end{equation*}
        is absolutely continuous, and
        \begin{equation*}
            \frac d{dt}\gamma^{-1}\circ\xi(t)=\frac{\dot\gamma\mleft(\gamma^{-1}\circ\xi(t)\mright)\dot\xi(t)}{|\dot\gamma(\gamma^{-1}\circ\xi(t))|^2},\qquad \text{for a.e. }t\in[0,T].
        \end{equation*}
    \end{lem}

    \begin{defin}
        Given an absolutely continuous curve $\zeta:[0,T']\to\Rds^N$, a curve $\xi:[0,T]\to\Rds^N$ is called a \emph{reparametrization} of $\zeta$ if there exists a nondecreasing surjective absolutely continuous function $\psi$ from $[0,T]$ onto $[0,T']$ with
        \begin{equation*}
            \xi(t)=\zeta\circ\psi(t),\qquad \text{for any }t\in[0,T].
        \end{equation*}
    \end{defin}

    Note that if $\xi$ is a reparametrization of $\zeta$, the converse property in general is not true for $\psi$ could have not strictly positive derivative for a.e.\ $t$, see Zarecki criterion for an absolutely continuous inverse in~\cite{Bernal18}. We have that reparametrizations are absolutely continuous:

    \begin{lem}\label{reparac}
        \emph{\cite[Corollary 4]{SerrinVarberg69}} Let $\zeta:[0,T']\to\Rds^N$ be a curve and $\psi:[0,T]\to[0,T']$ be absolutely continuous and nondecreasing. Then the reparametrization $\xi\equiv\zeta\circ\psi$ of $\zeta$ is absolutely continuous and
        \begin{equation*}
            \frac d{dt}\xi(t)=\dot\zeta(\psi(t))\dot\psi(t),\qquad \text{a.e.\ in }[0,T].
        \end{equation*}
    \end{lem}

    \begin{lem}\label{repsigma}
        If the curve $\xi:[0,T]\to\Gamma$ is a reparametrization of a curve $\zeta:[0,T']\to\Gamma$, then
        \begin{equation*}
            \int_0^{T'}\sigma_a\mleft(\xi,\dot\xi\mright)d\tau=\int_0^T\sigma_a\mleft(\zeta,\dot\zeta\mright)d\tau,\qquad \text{for every }a\in\Rds.
        \end{equation*}
    \end{lem}
    \begin{proof}
        It follows from the definition that $(x,q)\mapsto\sigma_a(x,q)$ is positively homogeneous on $q$, thus, if we let $\psi$ be the nondecreasing absolutely continuous function such that $\xi\equiv\zeta\circ\psi$ and consider the change of variable $r=\psi(\tau)$, we get from \cref{reparac} that, for every $a\in\Rds$,
        \begin{equation*}
            \int_0^T\sigma_a\mleft(\xi,\dot\xi\mright)d\tau=\int_0^T\sigma_a\mleft(\zeta\circ\psi,\dot\zeta\circ\psi\mright)\dot\psi(\tau)d\tau=\int_0^{T'}\sigma_a\mleft(\zeta,\dot\zeta\mright)dr.
        \end{equation*}
    \end{proof}

    The next \namecref{repconstspeed} comes from classical results of analysis in metric space, see \cite[Proposition 4.14]{Bernal18}.

    \begin{prop}\label{repconstspeed}
        Any curve defined on a bounded interval is a reparametrization of some curve $\zeta:[0,T]\to\Rds^N$ with constant speed, i.e., with $\mleft|\dot\zeta\mright|=\mathrm{constant}$ a.e..
    \end{prop}

    \begin{rem}\label{stpar}
        Let $\xi:[0,T]\to\Gamma$ be a simple curve. Since $\xi$ is absolutely continuous, there is a finite partition $\{t_0,\dotsc,t_m\}$ of the interval $[0,T]$ so that $t_0=0$, $t_m=T$ and
        \begin{equation*}
            \xi((t_{i-1},t_i))\subseteq\gamma_i((0,1)),\qquad \text{for each }i\in\{1,\dotsc,m\},
        \end{equation*}
        where $\gamma_i$ is an arc of $\Gamma$. \Cref{repconstspeed,invarccomp} yield that, for every $i\in\{1,\dotsc,m\}$, $\gamma^{-1}_i\circ\xi|_{[t_{i-1},t_i]}$ is the reparametrization of a curve $\eta_i$ with either $\dot\eta_i=1$ a.e.\ or $\dot\eta_i=-1$ a.e.. Possibly replacing $\gamma_i$ with $\wtgamma_i$, we can always assume that $\dot\eta_i=1$ a.e.. Setting
        \begin{equation*}
            \zeta=(\gamma_1\circ\eta_1)*\dotsb*(\gamma_m\circ\eta_m),
        \end{equation*}
        it is apparent that $\xi$ is a reparametrization of $\zeta$. This shows that any simple curve on $\Gamma$ is the reparametrization of a finite concatenation of sub-arcs, see \cref{subarcdef}.
    \end{rem}

    The \namecref{closedcurve} below is a simple consequence of~\cref{eq:critvaldef} and \cref{subsolchar}.

    \begin{lem}\label{closedcurve}
        Let $\xi:[0,T]\to\Gamma$ be a closed curve. Then, for any $a\ge c$,
        \begin{equation*}
            \int_0^T\sigma_a\mleft(\xi,\dot\xi\mright)d\tau\ge0.
        \end{equation*}
    \end{lem}

    \begin{prop}\label{lowcost}
        Given a nonconstant curve $\xi:[0,T]\to\Gamma$ there is a curve $\zeta:[0,T']\to\Gamma$, which is a finite concatenation of sub-arcs and has the same endpoints of $\xi$, so that
        \begin{equation*}
            \int_0^T\sigma_a\mleft(\xi,\dot\xi\mright)d\tau\ge\int_0^{T'}\sigma_a\mleft(\zeta,\dot\zeta\mright)d\tau,\qquad \text{whenever }a\ge c.
        \end{equation*}
        Furthermore, if $\xi$ is simple or non-closed then $\zeta$ is simple.
    \end{prop}
    \begin{proof}
        We preliminarily assume that $\xi$ is not closed. We define $E$ as the set made up of the nonempty intervals $[t_1,t_2)\subset[0,T]$ such that the restriction of $\xi$ to $[t_1,t_2]$ is a closed curve. These intervals are at most countable, therefore $E$ is a measurable set and
        \begin{equation*}
            \dot\xi_0(t)\coloneqq\mleft\{
            \begin{aligned}
                &\dot\xi(t),&& \text{if }t\in[0,T]\setminus E,\\
                &0,&& \text{if }t\in E,
            \end{aligned}
            \mright.
        \end{equation*}
        is a measurable function, i.e., is the derivative of an absolutely continuous curve $\xi_0:[0,T]\to\Gamma$. Thanks to \cref{repconstspeed} we have that $\xi_0$ is the reparametrization of a curve $\zeta:[0,T']\to\Gamma$ with constant speed. We point out that $\zeta$ is simple by construction, thus we can assume, in view of \cref{stpar}, that $\zeta$ is a concatenation of sub-arcs. Finally \cref{closedcurve,repsigma} yield
        \begin{equation*}
            \int_0^T\sigma_a\mleft(\xi,\dot\xi\mright)d\tau\ge\int_0^T\sigma_a\mleft(\xi_0,\dot\xi_0\mright)d\tau=\int_0^{T'}\sigma_a\mleft(\zeta,\dot\zeta\mright)d\tau,\qquad \text{whenever }a\ge c,
        \end{equation*}
        which proves our claim when $\xi$ is not closed. The case where $\xi$ is closed can be solved similarly, breaking $\xi$ into two non-closed curves.
    \end{proof}

    We can now prove the main results of this \lcnamecref{optcurve}.

    \mincurve*
    \begin{proof}
        The proof of~\ref{en:mincurve1} is trivial in view of \cref{closedcurve}, thus we focus on \cref{en:mincurve2}. We notice that by \cref{lowcost} the infimum in~\cref{eq:semidist} can be taken over the simple curves linking $y$ and $x$ which are a concatenation of sub-arcs. Since there is only a finite number of arcs, it is apparent that there is only a finite number of such curves. This concludes our proof.
    \end{proof}

    \sdiscont*
    \begin{proof}
        Fixed $a\ge c$ and $(y_1,x_1),(y_2,x_2)\in\Gamma^2$, let $\zeta:[0,T]\to\Gamma$ be an optimal curve for $S_a(y_2,x_2)$, which exists by \cref{mincurve}, $\xi_y:[0,T_y]\to\Gamma$ be a geodesic from $y_1$ to $y_2$ and $\xi_x:[0,T_x]\to\Gamma$ be a geodesic from $x_2$ to $x_1$. It is apparent that
        \begin{align*}
            S_a(y_1,x_1)-S_a(y_2,x_2)\le&\int_0^{T_y}\sigma_a\mleft(\xi_y,\dot\xi_y\mright)d\tau+\int_0^T\sigma_a\mleft(\zeta,\dot\zeta\mright)d\tau+\int_0^{T_x}\sigma_a\mleft(\xi_x,\dot\xi_x\mright)d\tau\\
            &-\int_0^T\sigma_a\mleft(\zeta,\dot\zeta\mright)d\tau\\
            \le&\int_0^{T_y}\sigma_a\mleft(\xi_y,\dot\xi_y\mright)d\tau+\int_0^{T_x}\sigma_a\mleft(\xi_x,\dot\xi_x\mright)d\tau.
        \end{align*}
        Setting
        \begin{equation*}
            \ell_a\coloneqq\max_{(x,q)\in T\Gamma,|q|\le1}\sigma_a(x,q),
        \end{equation*}
        we can then exploit the positive homogeneity of $q\mapsto\sigma_a(x,q)$ to get
        \begin{equation*}
            S_a(y_1,x_1)-S_a(y_2,x_2)\le\ell_a\int_0^{T_y}\mleft|\dot\xi_y(\tau)\mright|d\tau+\ell_a\int_0^{T_x}\mleft|\dot\xi_x(\tau)\mright|d\tau=\ell_a d_\Gamma(y_1,y_2)+\ell_a d_\Gamma(x_2,x_1).
        \end{equation*}
        Interchanging $(y_1,x_1)$, $(y_2,x_2)$ and since $(y_1,x_1)$, $(y_2,x_2)$ are arbitrary, this proves that $S_a$ is $\ell_a$--Lipschitz continuous.
    \end{proof}

    \printbibliography[heading=bibintoc]

\end{document}